\documentclass[11pt]{article}
\usepackage{bridges}
\usepackage{amsfonts,amssymb,amsthm,eucal,amsmath}
\usepackage{graphicx}
\usepackage{subfig}

\usepackage{wrapfig}
\usepackage{pinlabel, color}
\usepackage{multirow}

\usepackage{hyperref}

\newcommand{\R}{\mathbb{R}}

\newcommand{\Ha}{\mathbb{H}}
\newcommand{\SO}{\text{SO}}
\newcommand{\RP}{\mathbb{R}\text{P}}

\title{Hypernom: Mapping VR Headset Orientation to $S^3$}
\author{
\begin{tabular}{cc}
Vi Hart            & Andrea Hawksley \\
Communications Design Group    & Communications Design Group  \\
SAP Labs       & SAP Labs     \\
 & \\
Henry Segerman & Marc ten Bosch  \\
Department of Mathematics & MTB Design Works, Inc.\\
 Oklahoma State University\\
\end{tabular}}

\date{}

\begin{document}
\maketitle

\begin{abstract}
\emph{Hypernom} is a virtual reality game. The cells of a regular 4D polytope are radially projected to $S^3$, the sphere in 4D space, then stereographically projected to 3D space where they are viewed in the headset. The orientation of the headset is given by an element of the group $\SO(3)$, which is also a space that is double covered by $S^3$. In fact, the headset outputs a point of this double cover: a unit quaternion. The positions of the cells are multiplied by this quaternion before projection to 3D space, which moves the player through $S^3$. When the player is sufficiently close to a cell, they eat it. The aim of the game is to eat all of the cells of the polytope, which, roughly speaking, is achieved by moving one's head through all possible orientations, twice. 
 \end{abstract}

\section*{Introduction}

\emph{Hypernom} is a virtual reality game, available at \url{http://hypernom.com}. Although a virtual reality headset (such as an Oculus Rift) and game pad are required to appreciate the full experience, the game is also playable using an ordinary computer or smart phone. Instructions are provided at the link.

\begin{figure}
\centering

\subfloat[Cells of the 120-cell after stereographic projection, rendered to a virtual reality headset.]
{
\includegraphics[height=1.7in]{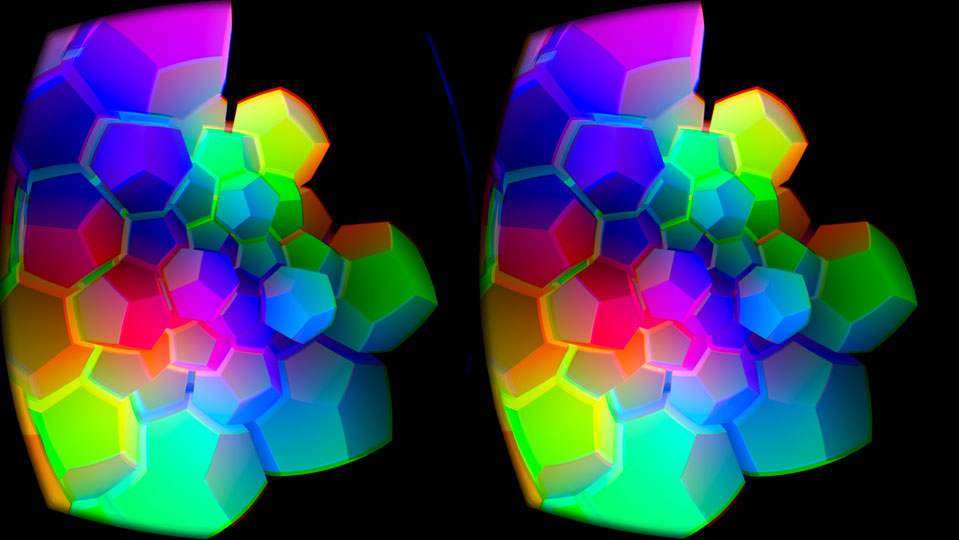}
\label{120cellVR2.jpg}
} 
\thinspace
\subfloat[Cells of the 600-cell after stereographic projection, rendered to a screen.]
{
\includegraphics[height=1.7in]{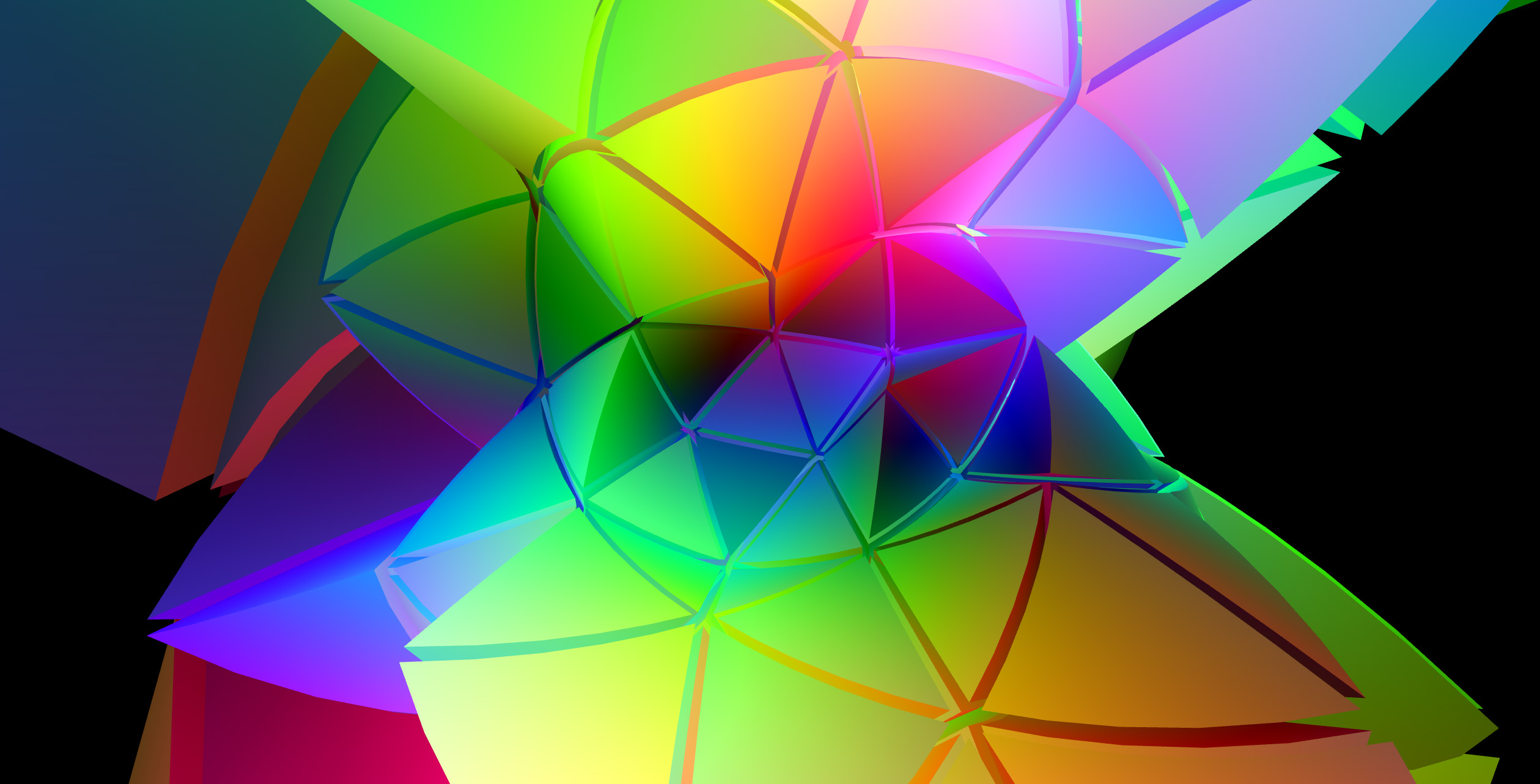}
\label{600cellRect.jpg}
} 

\caption{Screenshots from \emph{Hypernom}. The cells are coloured according to an algorithm based on the Hopf fibration, which is a map from $S^3$ to $S^2$, the sphere in $\R^3$. The sphere $S^2$ is then inscribed in the cube whose $x, y$, and $z$ coordinates are identified with brightness values for the red, green and blue components of a colour respectively. Thus, a point on the sphere gives a colour. We use this map from $S^3$ to $S^2$ to the colour cube in two ways. First, each point of a cell is at some position in $S^3$, and this gives a ``base colour'' for that point of the cell. Second, a point on the boundary surface of a cell has a normal vector, which is a unit vector in $\R^4$, and therefore is also a point in $S^3$. This point gives the ``shading'' effect -- it is added to the base colour to give the final result.}
\label{screenshots}
\end{figure}

The game is played on a choice of one of the six regular four-dimensional polytopes.  These are the 5-cell, the hypercube or 8-cell, the 16-cell, the 24-cell, the 120-cell and the 600-cell (see \cite[page~136]{Coxeter1973}). As in previous work of the third author \cite{bridges2012:103}, we view the cells of these polytopes in $\R^3$ by first radially projecting the polytope in $\R^4$ onto the unit sphere, $S^3=\{(w,x,y,z)\mid w^2+x^2+y^2+z^2=1\}$, then stereographically projecting into $\R^3$. The result is rendered on the screen or virtual reality headset using standard 3D graphics techniques, with the camera positioned at the origin of $\R^3$ by default. See Figure \ref{screenshots}.

The main innovation in \emph{Hypernom} is the unusual way in which we interpret the orientation information that the virtual reality headset feeds back to the computer. Traditionally, headset orientation information is used to control the orientation of a virtual camera in an artificial scene, to simulate as accurately as possible what the viewer would see when looking around the virtual environment.
We rotate a virtual camera in this way, but also use the orientation data to define the player's position in space, as described in the next section.

When a cell of the polytope comes sufficiently close to the virtual camera (as measured by distance in $S^3$, before stereographic projection), the player ``eats'' it -- the cell disappears from view, revealing any cells behind it. The aim of the game is to eat all of the cells of the polytope.

\section*{Moving through $S^3$}\label{moving_through_S3}

The \emph{real quaternions} $\Ha$ are a number system that extends the real numbers, similarly to the way in which the complex numbers also extend the real numbers.
$$\Ha = \{w +xi+yj+zk \mid w,x,y,z\in\R, i^2=j^2=k^2=ijk=-1 \}$$
The \emph{unit} quaternions are those that, viewed as points of $\R^4$, are at distance one from the origin. These are then in one-to-one correspondence with the points of $S^3$. Said another way, we can give the points of $S^3$ the structure of a multiplicative group by viewing them as unit quaternions.

 \begin{wrapfigure}[12]{l}{0.4\textwidth}
 \centering
 \vspace{-5pt}
 \labellist
 \small\hair 2pt
 \pinlabel roll at 1800 707 
 \pinlabel yaw at 1680 370
 \pinlabel pitch at 1130 100 
 \endlabellist
 \includegraphics[width=0.3\textwidth]{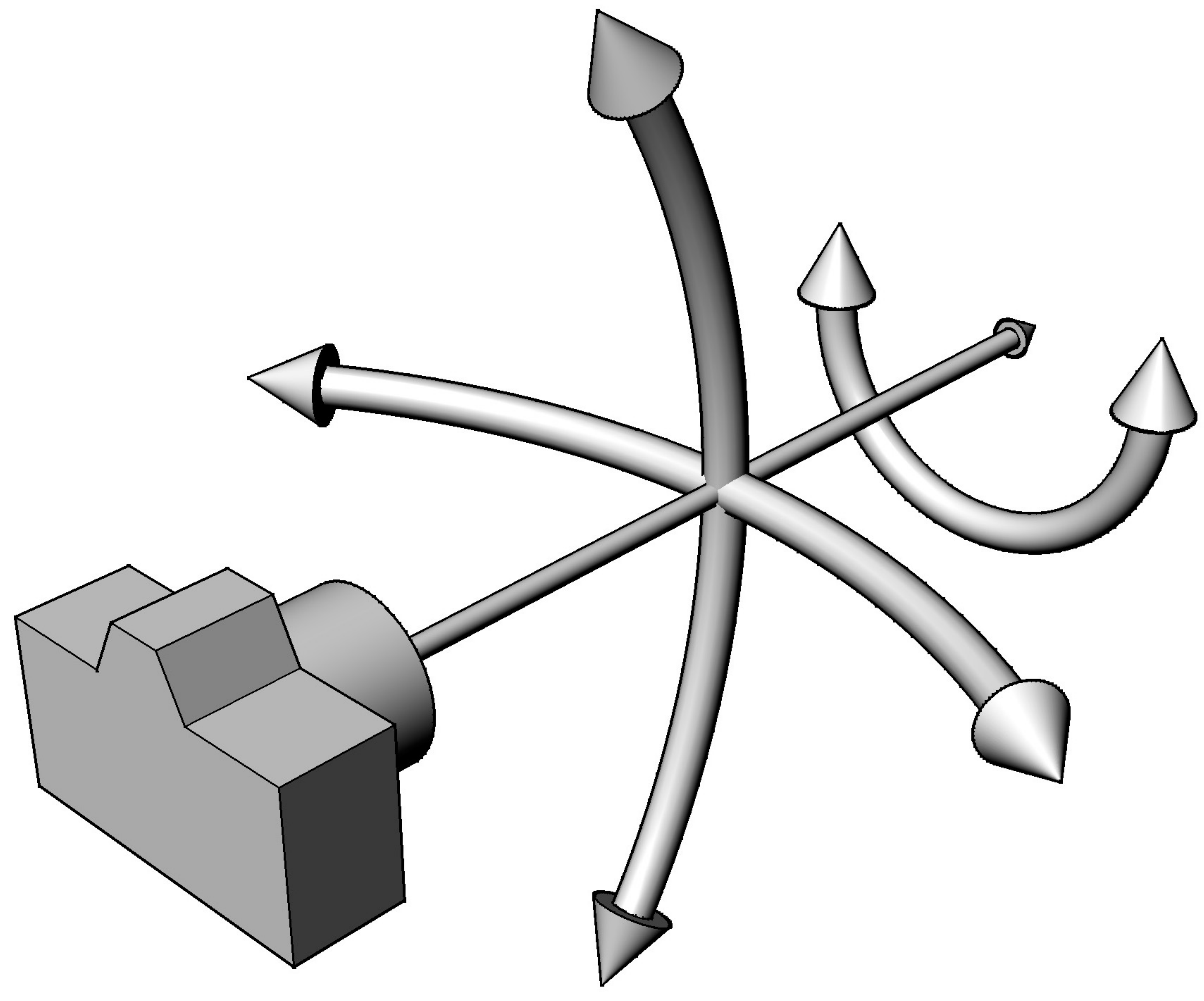}
 \caption{Changes of orientation relative to the current orientation.} 
 \label{pitch_yaw_roll}
 \end{wrapfigure}
 
The set of possible orientations of an object (for example, your head, or a camera) in $\R^3$ is a three-dimensional space. There are always three degrees of freedom in orientation: one can change the pitch, yaw and roll. See Figure \ref{pitch_yaw_roll}. This set of possible orientations is $\SO(3)$, the group of three by three orthogonal matrices with determinant one. As a topological space, $\SO(3)$ is homeomorphic to $\RP^3$, the space obtained from $S^3$ by identifying antipodal points. To see this somewhat intuitively, to each point inside of the 3-dimensional ball of radius $\pi$ centered at the origin, we associate a change of orientation as follows (see also \cite[Exercise 2.7.7]{thurston_book}). The origin, $0$, corresponds to no change in orientation. For points $p \neq 0$, we view  $p$ as a vector pointing outwards from $0$, which gives us an axis we rotate around by an angle given by the length of $p$ (with some fixed choice of direction of rotation). By Euler's rotation theorem~\cite{palais_palais_rodi}, every possible change of orientation is given within this ball. Antipodal points on the boundary of the ball give the same change in orientation, since these points correspond to rotation around the same axis by $\pi$ and $-\pi$. So, the space of possible orientations is formed from the ball by gluing opposite points on its boundary, giving $\RP^3$. This is analogous to gluing opposite points on the boundary of a disk to obtain $\RP^2$, the space obtained from $S^2$ by identifying antipodal points.

In modern computer graphics, orientations in 3D space are often encoded using unit quaternions (i.e., elements of $S^3$) rather than matrices in $\SO(3)$ or other representations. The double covering of $\SO(3)$ by $S^3$ (see \cite[page~106]{montesinos}) is a minor issue: what this means is that an orientation can be encoded using either of two quaternions which are negatives of each other. The orientation of a virtual reality headset is sensed using gyroscopes and accelerometers, then this data is converted into the quaternion representation. The choice of which of the two quaternions to use is determined by assuming continuity through time -- the quaternion we are sensing now should be close to the quaternion we sensed last frame.

So, the headset supplies us with a unit quaternion, and the points in our space $S^3$ are also unit quaternions. We act on $S^3$ (and so the cells of our chosen polytope) by multiplying all of its points by the headset quaternion, before stereographically projecting to $\R^3$. So, as the viewer rotates his or her head around, the cells of the polytope move and twist. Although the virtual camera is always positioned at the origin of $\R^3$, the visual effect is of the viewer moving through the space $S^3$.

\section*{Discussion}

\begin{figure}
\centering

\subfloat[]
{
\includegraphics[width=0.45\textwidth]{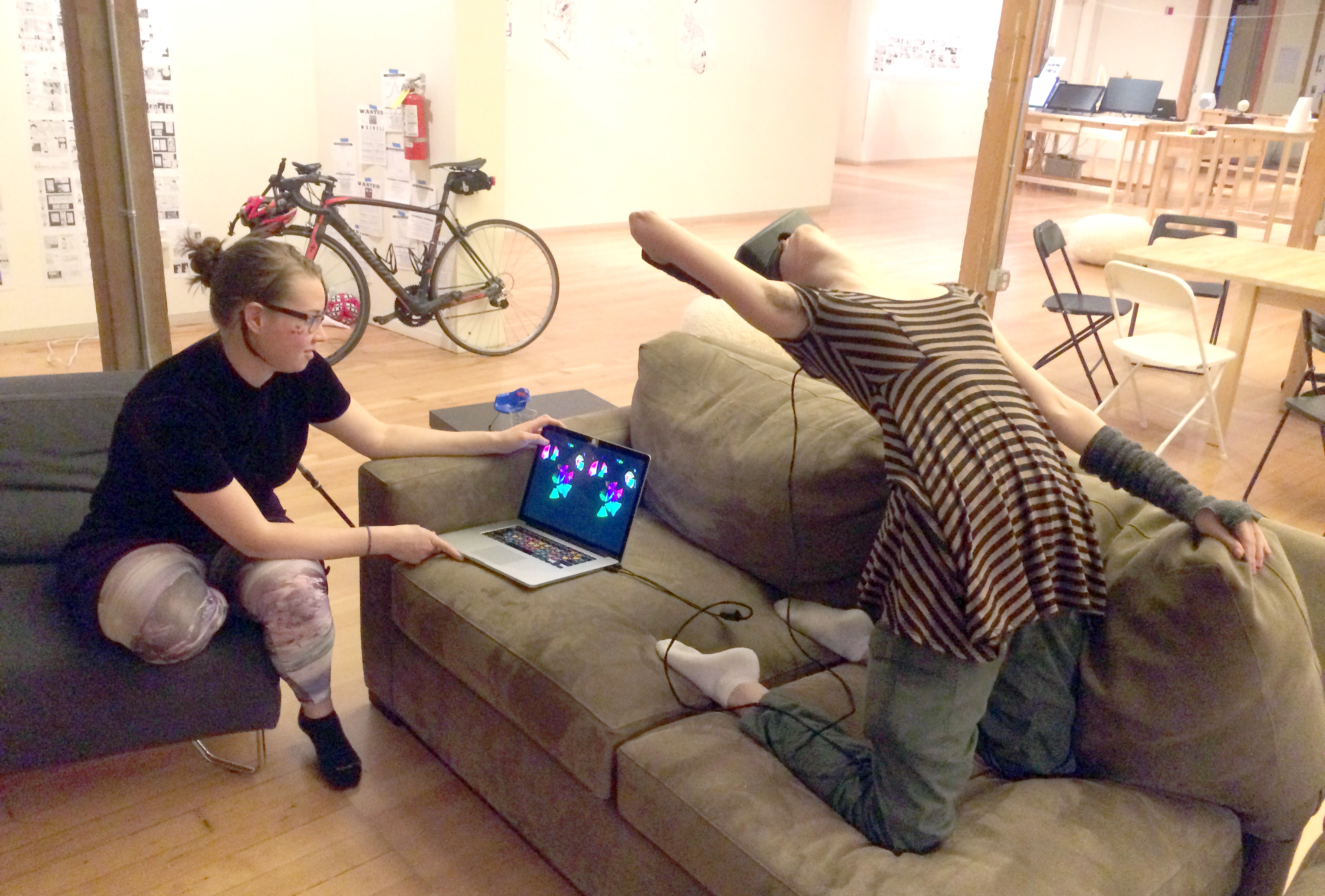}
\label{vi_couch1}
} 
\qquad
\subfloat[]
{
\includegraphics[width=0.45\textwidth]{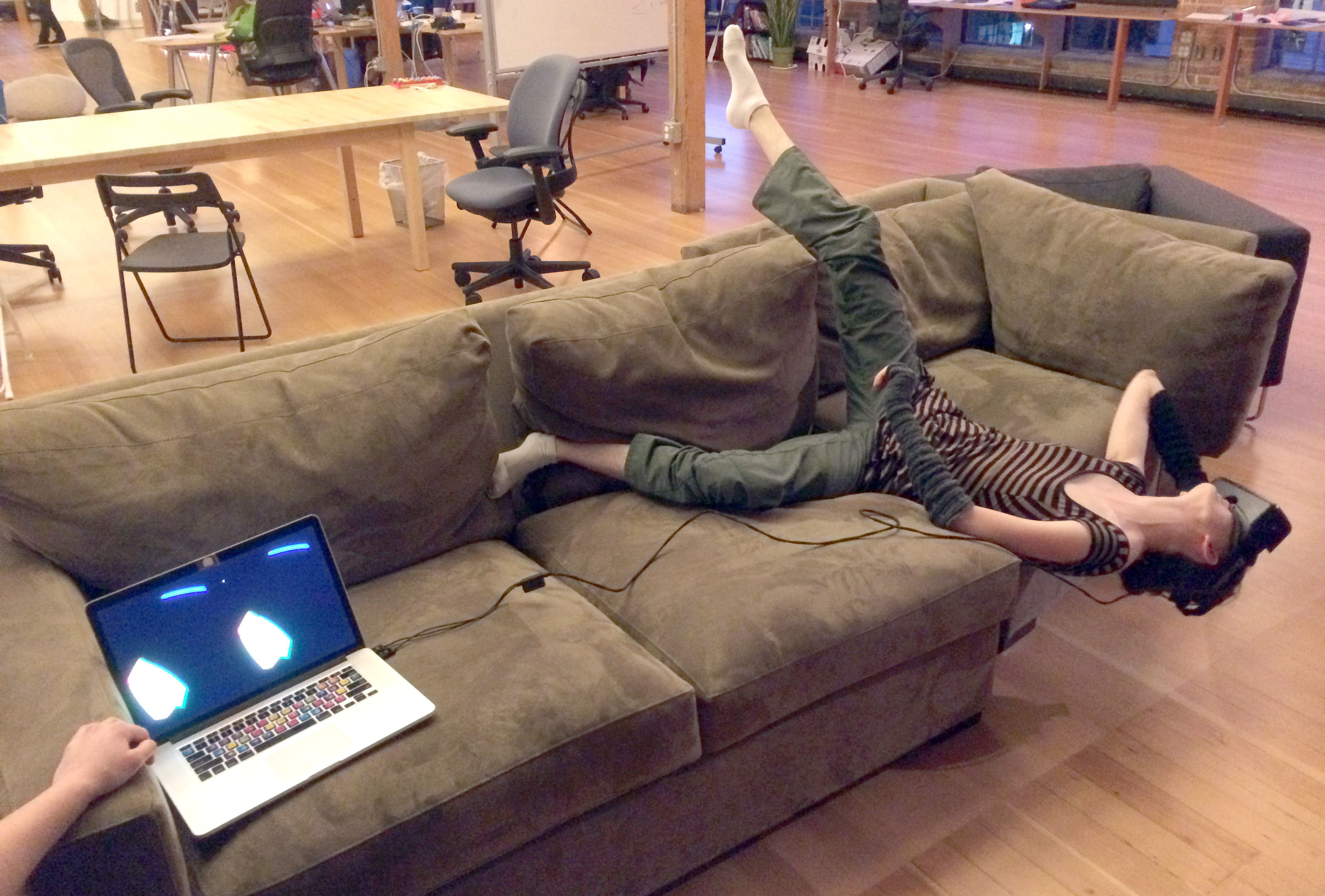}
\label{vi_couch2}
} 
\caption{Winning \emph{Hypernom} requires moving one's head close to every possible orientation.}
\label{vi_couch}
\end{figure}

How does one win the game? In one sense, the answer is relatively simple: each of the regular polytopes is a regular subdivision of $S^3$. To eat a cell, one must move sufficiently close to the center of each cell, the distance required depending on the size of the cells. Particularly for the polytopes with a large number of cells, this means that we have to get close to every point of $S^3$, our path approximating an $S^3$-filling curve, for example following a Hamiltonian path on the edges of the dual polytope~\cite{sequin_hamiltonian}. Since the player controls their ``position'' using the orientation of their head, they must move their head near to every possible orientation. In fact they must do this twice, in order to account for the double covering of $\SO(3)$ by $S^3$. This is not very easy to do for the 120- and 600-cells -- although many cells can be eaten by looking around or spinning on the spot, there are many others that require the player's head to be upside down. We recommend playing the game on a comfortable couch which may be flopped around on. See Figure \ref{vi_couch}.

The fact that $S^3$ is the double cover of $\SO(3)$ can be handily demonstrated within \emph{Hypernom}. If the player turns around on the spot by 360 degrees, they travel through the polytope eating cells, but do not get back to where they started. However, making another full turn, for a total rotation of 720 degrees does return the player to their starting position.

When players first try the game, they quickly pick up a sense of what kind of motions of their head perform what kind of movement through the $S^3$. From a position looking forward while standing up, changing your orientation by altering yaw also moves you vertically; altering pitch also moves you horizontally; altering roll also moves you forwards or backwards. In fact, the correspondence between the change of head orientation relative to the current facing and the perceived motion of cells in $S^3$, holds no matter which way one is looking.  So, one strategy is to spot a cell in front of you and attempt to roll one's head to the right while keeping the direction of your facing fixed. This brings the cell towards you so you can eat it. This is often easier said than done however, potentially requiring handstands, cartwheels and the like.

\section*{Motivation and Artistic Choices}

It may seem somewhat arbitrary to design a VR game that uses headset orientation data as a quaternion to map your position in $S^3$ to eat cells of regular polychora, but given the context of mathematical art and VR research the authors are immersed in, it seemed obvious, even necessary, to design exactly this game.

Regular polytopes in 4D, like the 3D platonic solids, are mathematical objects that have motivated many people to visualize and share their beauty. Visualizing 4D objects in 3D has been done in many ways, each with strengths and weaknesses. 4D objects rendered and projected onto a computer screen can be manipulated in ways physical sculpture cannot. However, they lack the physical presence and 3-dimensionality that a physical sculpture or 3D print can provide. This becomes all the more important when working a dimension short. The authors are intrigued by the potential of virtual reality to combine the mutability of rendered digital objects with the sense of physical space that previously only physical objects could provide.

Projects leading up to \emph{Hypernom} immersed the user in hyperbolic and spherical spaces. Regular tilings are helpful aids for understanding a curved space: while we cannot necessarily trust our eyes to see things as they really are in a non-Euclidean space, we can trust our understanding of a single tile and know that all other tiles are identical, regardless of how warped they appear. The decision to fill $S^3$ with regular polychora in this game came after the decision to use headset orientation to navigate $S^3$. Being familiar with the relationship between $\SO(3)$ and $S^3$, and having used both independently in previous VR work \cite{monkeys}, we were curious whether applying one to the other would result in something that made sense to the human brain and led to further understanding of the space, rather than confusion. The mapping is surprisingly comprehensible, with head motions feeling like a method of navigation through a space. 

Collecting or eating objects is a classic game mechanic for showing progress, from PacMan to Mario Bros. In this game, it serves a dual purpose: it shows you how much of $S^3$ you've covered (to some degree of accuracy), and, as you're embedded in the space itself, eating nearby cells clears the way to see more of the space. A mechanic where you merely highlight, or even create, cells, might work if we were 4D creatures looking in at $S^3$. But being merely 3D, a subtractive method is the obvious choice. Given the history of the game mechanic, the design decision was clear.
The different regular polychora give us different levels of granularity, in terms of how much of $S^3$ the player needs to cover to win -- this allows us to vary the difficulty from very low accuracy for the 5-cell to fairly high accuracy for the 600-cell. Thus, even the choice of including all regular polychora made sense for both aesthetic and gameplay reasons.

\bibliographystyle{plain} 
\bibliography{hypernom_biblio}

\end{document}